\documentclass{amsart}
\usepackage{amsmath}
\usepackage{amsfonts}

\newtheorem{theorem}{Theorem}
\theoremstyle{plain}

\newtheorem{corollary}{Corollary}

\newtheorem{proposition}{Proposition}

\numberwithin{equation}{section}
\input{tcilatex}

\begin{document}
\title[On the $q$-Analogue of $p$-Adic $\log $ gamma type functions ]{$q$%
-Analogue of $p$-Adic $\log $ $\Gamma $ type functions associated with
Modified $q$-extension of Genocchi numbers with weight $\alpha $ and $\beta $
}
\author{Serkan Arac\i }
\address{University of Gaziantep, Faculty of Science and Arts, Department of
Mathematics, 27310 Gaziantep, TURKEY}
\email{mtsrkn@hotmail.com}
\author{Mehmet Acikgoz}
\address{University of Gaziantep, Faculty of Science and Arts, Department of
Mathematics, 27310 Gaziantep, TURKEY}
\email{acikgoz@gantep.edu.tr}
\date{January 20, 2012}

\begin{abstract}
The fundamental aim of this paper is to describe $q$-Analogue of $p$-adic $%
\log $ gamma functions with weight alpha and beta. Moreover, we give
relationship between $p$-adic $q$-$\log $ gamma funtions with weight ($%
\alpha ,\beta $) and $q$-extension of Genocchi numbers with weight alpha and
beta and modified $q$-Euler numbers with weight $\alpha $
\end{abstract}

\maketitle

\section{Introduction}

Assume that $p$ be a fixed odd prime number. Throughout this paper $%
\mathbb{Z}
,$ $%
\mathbb{Z}
_{p},$ $%
\mathbb{Q}
_{p}$ and $%
\mathbb{C}
_{p}$ will denote by the ring of integers, the field of $p$-adic rational
numbers and the completion of the algebraic closure of $%
\mathbb{Q}
_{p},$ respectively. Also we denote $%
\mathbb{N}
^{\ast }=%
\mathbb{N}
\cup \left\{ 0\right\} $ and $\exp \left( x\right) =e^{x}.$ Let $v_{p}:%
\mathbb{C}
_{p}\rightarrow 
\mathbb{Q}
\cup \left\{ \infty \right\} \left( 
\mathbb{Q}
\text{ is the field of rational numbers}\right) $ denote the $p$-adic
valuation of $%
\mathbb{C}
_{p}$ normalized so that $v_{p}\left( p\right) =1$. The absolute value on $%
\mathbb{C}
_{p}$ will be denoted as $\left\vert .\right\vert _{p}$, and $\left\vert
x\right\vert _{p}=p^{-v_{p}\left( x\right) }$ for $x\in 
\mathbb{C}
_{p}.$ When one talks of $q$-extensions, $q$ is considered in many ways,
e.g. as an indeterminate, a complex number $q\in 
\mathbb{C}
,$ or a $p$-adic number $q\in 
\mathbb{C}
_{p},$ If $q\in 
\mathbb{C}
$ we assume that $\left\vert q\right\vert <1.$ If $q\in 
\mathbb{C}
_{p},$ we assume $\left\vert 1-q\right\vert _{p}<p^{-\frac{1}{p-1}},$ so
that $q^{x}=\exp \left( x\log q\right) $ for $\left\vert x\right\vert
_{p}\leq 1.$ We use the following notation 
\begin{equation}
\left[ x\right] _{q}=\frac{1-q^{x}}{1-q},\text{ \ }\left[ x\right] _{-q}=%
\frac{1-\left( -q\right) ^{x}}{1+q}  \label{equation 1}
\end{equation}

where $\lim_{q\rightarrow 1}\left[ x\right] _{q}=x;$ cf. [1-24].

For a fixed positive integer $d$ with $\left( d,f\right) =1,$ we set%
\begin{eqnarray*}
X &=&X_{d}=\lim_{\overleftarrow{N}}%
\mathbb{Z}
/dp^{N}%
\mathbb{Z}
, \\
X^{\ast } &=&\underset{\underset{\left( a,p\right) =1}{0<a<dp}}{\cup }a+dp%
\mathbb{Z}
_{p}
\end{eqnarray*}

and%
\begin{equation*}
a+dp^{N}%
\mathbb{Z}
_{p}=\left\{ x\in X\mid x\equiv a\left( \func{mod}dp^{N}\right) \right\} ,
\end{equation*}

where $a\in 
\mathbb{Z}
$ satisfies the condition $0\leq a<dp^{N}.$

It is known that 
\begin{equation*}
\mu _{q}\left( x+p^{N}%
\mathbb{Z}
_{p}\right) =\frac{q^{x}}{\left[ p^{N}\right] _{q}}
\end{equation*}

is a distribution on $X$ for $q\in 
\mathbb{C}
_{p}$ with $\left\vert 1-q\right\vert _{p}\leq 1.$

Let $UD\left( 
\mathbb{Z}
_{p}\right) $ be the set of uniformly differentiable function on $%
\mathbb{Z}
_{p}.$ We say that $f$ is a uniformly differentiable function at a point $%
a\in 
\mathbb{Z}
_{p},$ if the difference quotient 
\begin{equation*}
F_{f}\left( x,y\right) =\frac{f\left( x\right) -f\left( y\right) }{x-y}
\end{equation*}

has a limit $f^{%
{\acute{}}%
}\left( a\right) $ as $\left( x,y\right) \rightarrow \left( a,a\right) $ and
denote this by $f\in UD\left( 
\mathbb{Z}
_{p}\right) .$ The $p$-adic $q$-integral of the function $f\in UD\left( 
\mathbb{Z}
_{p}\right) $ is defined by%
\begin{equation}
I_{q}\left( f\right) =\int_{%
\mathbb{Z}
_{p}}f\left( x\right) d\mu _{q}\left( x\right) =\lim_{N\rightarrow \infty }%
\frac{1}{\left[ p^{N}\right] _{q}}\sum_{x=0}^{p^{N}-1}f\left( x\right) q^{x}
\label{equation 2}
\end{equation}

The bosonic integral is considered by Kim as the bosonic limit $q\rightarrow
1,$ $I_{1}\left( f\right) =\lim_{q\rightarrow 1}I_{q}\left( f\right) .$
Similarly, the $p$-adic fermionic integration on $%
\mathbb{Z}
_{p}$ \ defined by Kim as follows:%
\begin{equation*}
I_{-q}\left( f\right) =\lim_{q\rightarrow -q}I_{q}\left( f\right) =\int_{%
\mathbb{Z}
_{p}}f\left( x\right) d\mu _{-q}\left( x\right)
\end{equation*}

Let $q\rightarrow 1,$ then we have $p$-adic fermionic integral on $%
\mathbb{Z}
_{p}$ as follows:%
\begin{equation*}
I_{-1}\left( f\right) =\lim_{q\rightarrow -1}I_{q}\left( f\right)
=\lim_{N\rightarrow \infty }\sum_{x=0}^{p^{N}-1}f\left( x\right) \left(
-1\right) ^{x}.
\end{equation*}

Stirling asymptotic series are defined by%
\begin{equation}
\log \left( \frac{\Gamma \left( x+1\right) }{\sqrt{2\pi }}\right) =\left( x-%
\frac{1}{2}\right) \log x+\sum_{n=1}^{\infty }\frac{\left( -1\right) ^{n+1}}{%
n\left( n+1\right) }\frac{B_{n+1}}{x^{n}}-x  \label{equation 11}
\end{equation}

where $B_{n}$ are familiar $n$-th Bernoulli numbers cf. [6, 8, 9, 25].

Recently, Araci et al. defined modified $q$-Genocchi numbers and polynomials
with weight $\alpha $ and $\beta $ in [4, 5] by the means of generating
function:%
\begin{equation}
\sum_{n=0}^{\infty }g_{n,q}^{\left( \alpha ,\beta \right) }\left( x\right) 
\frac{t^{n}}{n!}=t\int_{%
\mathbb{Z}
_{p}}q^{-\beta \xi }e^{\left[ x+\xi \right] _{q^{\alpha }}t}d\mu _{-q^{\beta
}}\left( \xi \right)  \label{equation 3}
\end{equation}

So from above, we easily get Witt's formula of modified$\ q$-Genocchi
numbers and polynomials with weight $\alpha $ and $\beta $ as follows:%
\begin{equation}
\frac{g_{n+1,q}^{\left( \alpha ,\beta \right) }\left( x\right) }{n+1}=\int_{%
\mathbb{Z}
_{p}}q^{-\beta \xi }\left[ x+\xi \right] _{q^{\alpha }}^{n}d\mu _{-q^{\beta
}}\left( \xi \right)  \label{equation 12}
\end{equation}

where $g_{n,q}^{\left( \alpha ,\beta \right) }\left( 0\right)
:=g_{n,q}^{\left( \alpha ,\beta \right) }$ are modified $q$ extension of
Genocchi numbers with weight $\alpha $ and $\beta $ cf. [4,5].

In \cite{Rim}, Rim and Jeong are defined modified $q$-Euler numbers with
weight $\alpha $ as follows:%
\begin{equation}
\widetilde{\xi }_{n,q}^{\left( \alpha \right) }=\int_{%
\mathbb{Z}
_{p}}q^{-t}\left[ t\right] _{q^{\alpha }}d\mu _{-q}\left( t\right)
\label{equation 24}
\end{equation}

\bigskip From expressions of (\ref{equation 12}) and (\ref{equation 24}), we
get the following Proposition 1:

\begin{proposition}
The following%
\begin{equation}
\widetilde{\xi }_{n,q}^{\left( \alpha \right) }=\frac{g_{n+1,q}^{\left(
\alpha ,1\right) }}{n+1}  \label{equation 25}
\end{equation}%
is true.
\end{proposition}

In previous paper \cite{Araci 6}, Araci, Acikgoz and Park introduced
weighted $q$-Analogue of $p$-Adic $\log $ gamma type functions and they
derived some interesting identities in Analytic Numbers Theory and in $p$%
-Adic Analysis. They were motivated from paper of T. Kim by "\textit{On a }$%
q $\textit{-analogue of the }$p$\textit{-adic log gamma functions and
related integrals}, \textit{J. Number Theory}, \textit{76 (1999), no. 2,
320-329}." \ We also introduce $q$-Analogue of $p$-Adic $\log $ gamma type
function with weight $\alpha $ and $\beta .$ We derive in this paper some
interesting identities this type of functions.

\begin{center}
\textbf{On} \textbf{p-adic }$\log $ $\Gamma $ \textbf{function with weight }$%
\alpha $ and $\beta $
\end{center}

In this part, from (\ref{equation 2}), we begin with the following nice
identity:%
\begin{equation}
I_{-q}^{\left( \beta \right) }\left( q^{-\beta x}f_{n}\right) +\left(
-1\right) ^{n-1}I_{-q}^{\left( \beta \right) }\left( q^{-\beta x}f\right) = 
\left[ 2\right] _{q^{\beta }}\sum_{l=0}^{n-1}\left( -1\right)
^{n-1-l}f\left( l\right)  \label{equation 6}
\end{equation}

where $f_{n}\left( x\right) =f\left( x+n\right) $ and $n\in 
\mathbb{N}
$ (see \cite{Araci 4}).

In particular for $n=1$ into (\ref{equation 6}), we easily see that%
\begin{equation}
I_{-q}^{\left( \beta \right) }\left( q^{-\beta x}f_{1}\right)
+I_{-q}^{\left( \beta \right) }\left( q^{-\beta x}f\right) =\left[ 2\right]
_{q^{\beta }}f\left( 0\right) .  \label{equation 7}
\end{equation}

With the simple application, it is easy to indicate as follows:%
\begin{equation}
\left( \left( 1+x\right) \log \left( 1+x\right) \right) ^{%
{\acute{}}%
}=1+\log \left( 1+x\right) =1+\sum_{n=1}^{\infty }\frac{\left( -1\right)
^{n+1}}{n\left( n+1\right) }x^{n}  \label{equation 15}
\end{equation}

where $\left( \left( 1+x\right) \log \left( 1+x\right) \right) ^{%
{\acute{}}%
}=\frac{d}{dx}\left( \left( 1+x\right) \log \left( 1+x\right) \right) $

By expression of (\ref{equation 15}), we can derive%
\begin{equation}
\left( 1+x\right) \log \left( 1+x\right) =\sum_{n=1}^{\infty }\frac{\left(
-1\right) ^{n+1}}{n\left( n+1\right) }x^{n+1}+x+c,\text{ where }c\text{ is
constant.}  \label{equation 16}
\end{equation}

If we take $x=0,$ so we get $c=0.$ By expression of (\ref{equation 15}) and (%
\ref{equation 16}), we easily see that,%
\begin{equation}
\left( 1+x\right) \log \left( 1+x\right) =\sum_{n=1}^{\infty }\frac{\left(
-1\right) ^{n+1}}{n\left( n+1\right) }x^{n+1}+x.  \label{equation 17}
\end{equation}

It is considered by T. Kim for $q$-analogue of $p$ adic locally analytic
function on $%
\mathbb{C}
_{p}\backslash 
\mathbb{Z}
_{p}$ as follows:%
\begin{equation}
G_{p,q}\left( x\right) =\int_{%
\mathbb{Z}
_{p}}\left[ x+\xi \right] _{q}\left( \log \left[ x+\xi \right] _{q}-1\right)
d\mu _{-q}\left( \xi \right) \text{ (for detail, see[5,6]).}
\label{equation 18}
\end{equation}

By the same motivation of (\ref{equation 18}), in previous paper \cite{Araci
6}, $q$-analogue of $p$-adic locally analytic function on $%
\mathbb{C}
_{p}\backslash 
\mathbb{Z}
_{p}$ with weight $\alpha $\ is considered 
\begin{equation}
G_{p,q}^{\left( \alpha \right) }\left( x\right) =\int_{%
\mathbb{Z}
_{p}}\left[ x+\xi \right] _{q^{\alpha }}\left( \log \left[ x+\xi \right]
_{q^{\alpha }}-1\right) d\mu _{-q}\left( \xi \right) \text{ }
\label{equation 19}
\end{equation}

In particular $\alpha =1$ into (\ref{equation 19}), we easily see that, $%
G_{p,q}^{\left( 1\right) }\left( x\right) =G_{p,q}\left( x\right) .$

With the same manner, we introduce $q$-Analoge of $p$-adic locally analytic
function on $%
\mathbb{C}
_{p}\backslash 
\mathbb{Z}
_{p}$ with weight $\alpha $ and $\beta $ as follows:%
\begin{equation}
G_{p,q}^{\left( \alpha ,\beta \right) }\left( x\right) =\int_{%
\mathbb{Z}
_{p}}q^{-\beta \xi }\left[ x+\xi \right] _{q^{\alpha }}\left( \log \left[
x+\xi \right] _{q^{\alpha }}-1\right) d\mu _{-q^{\beta }}\left( \xi \right) 
\text{ }  \label{equation 22}
\end{equation}

From expressions of (\ref{equation 7}) and\ (\ref{equation 20}), we state
the following Theorem:

\begin{theorem}
The following identity holds:%
\begin{equation*}
G_{p,q}^{\left( \alpha ,\beta \right) }\left( x+1\right) +G_{p,q}^{\left(
\alpha ,\beta \right) }\left( x\right) =\left[ 2\right] _{q^{\beta }}\left[ x%
\right] _{q^{\alpha }}\left( \log \left[ x\right] _{q^{\alpha }}-1\right) .
\end{equation*}%
\ 
\end{theorem}

It is easy to show that,%
\begin{eqnarray}
\left[ x+\xi \right] _{q^{\alpha }} &=&\frac{1-q^{\alpha \left( x+\xi
\right) }}{1-q^{\alpha }}  \label{equation 20} \\
&=&\frac{1-q^{\alpha x}+q^{\alpha x}-q^{\alpha \left( x+\xi \right) }}{%
1-q^{\alpha }}  \notag \\
&=&\left( \frac{1-q^{\alpha x}}{1-q^{\alpha }}\right) +q^{\alpha x}\left( 
\frac{1-q^{\alpha \xi }}{1-q^{\alpha }}\right)  \notag \\
&=&\left[ x\right] _{q^{\alpha }}+q^{\alpha x}\left[ \xi \right] _{q^{\alpha
}}  \notag
\end{eqnarray}

Substituting $x\rightarrow \frac{q^{\alpha x}\left[ \xi \right] _{q^{\alpha
}}}{\left[ x\right] _{q^{\alpha }}}$ into (\ref{equation 17}) and by using (%
\ref{equation 20}), we get interesting formula:%
\begin{equation}
\left[ x+\xi \right] _{q^{\alpha }}\left( \log \left[ x+\xi \right]
_{q^{\alpha }}-1\right) =\left( \left[ x\right] _{q^{\alpha }}+q^{\alpha x}%
\left[ \xi \right] _{q^{\alpha }}\right) \log \left[ x\right] _{q^{\alpha
}}+\sum_{n=1}^{\infty }\frac{\left( -q^{\alpha x}\right) ^{n+1}}{n(n+1)}%
\frac{\left[ \xi \right] _{q^{\alpha }}^{n+1}}{\left[ x\right] _{q^{\alpha
}}^{n}}-\left[ x\right] _{q^{\alpha }}  \label{equation 21}
\end{equation}

If we substitute $\alpha =1$ into (\ref{equation 21}), we get Kim's $q$%
-Analogue of $p$-adic $\log $ gamma fuction (for detail, see[8]).

From expression of (\ref{equation 2}) and (\ref{equation 21}), we obtain
worthwhile and interesting theorems as follows:

\begin{theorem}
For $x\in 
\mathbb{C}
_{p}\backslash 
\mathbb{Z}
_{p}$ the following 
\begin{equation}
G_{p,q}^{\left( \alpha ,\beta \right) }\left( x\right) =\left( \frac{\left[ 2%
\right] _{q^{\beta }}}{2}\left[ x\right] _{q^{\alpha }}+q^{\alpha x}\frac{%
g_{2,q}^{\left( \alpha ,\beta \right) }}{2}\right) \log \left[ x\right]
_{q^{\alpha }}+\sum_{n=1}^{\infty }\frac{\left( -q^{\alpha x}\right) ^{n+1}}{%
n\left( n+1\right) \left( n+2\right) }\frac{g_{n+1,q}^{\left( \alpha ,\beta
\right) }}{\left[ x\right] _{q^{\alpha }}^{n}}-\left[ x\right] _{q^{\alpha }}%
\frac{\left[ 2\right] _{q^{\beta }}}{2}  \label{equation 26}
\end{equation}%
is true.\ 
\end{theorem}

\begin{corollary}
Taking $q\rightarrow 1$ into (\ref{equation 26}), we get nice identity:%
\begin{equation*}
G_{p,1}^{\left( \alpha ,\beta \right) }\left( x\right) =\left( x+\frac{G_{2}%
}{2}\right) \log x+\sum_{n=1}^{\infty }\frac{\left( -1\right) ^{n+1}}{%
n\left( n+1\right) \left( n+2\right) }\frac{G_{n+1}}{x}-x
\end{equation*}%
where $G_{n}$ are called famous Genocchi numbers.
\end{corollary}

\begin{theorem}
The following nice identity%
\begin{equation}
G_{p,q}^{\left( \alpha ,1\right) }\left( x\right) =\left( \frac{\left[ 2%
\right] _{q}}{2}\left[ x\right] _{q^{\alpha }}+q^{\alpha x}\widetilde{\xi }%
_{1,q}^{\left( \alpha \right) }\right) \log \left[ x\right] _{q^{\alpha
}}+\sum_{n=1}^{\infty }\frac{\left( -q^{\alpha x}\right) ^{n+1}}{n\left(
n+1\right) }\frac{\widetilde{\xi }_{n,q}^{\left( \alpha \right) }}{\left[ x%
\right] _{q^{\alpha }}^{n}}-\frac{\left[ 2\right] _{q}}{2}\left[ x\right]
_{q^{\alpha }}  \label{equation 27}
\end{equation}%
is true.
\end{theorem}

\begin{corollary}
Putting $q\rightarrow 1$ into (\ref{equation 27}), we have the following
identity:%
\begin{equation*}
G_{p,1}^{\left( \alpha ,\beta \right) }\left( x\right) =\left(
x+E_{1}\right) \log x+\sum_{n=1}^{\infty }\frac{\left( -1\right) ^{n+1}}{%
n\left( n+1\right) }\frac{E_{n}}{x^{n}}-x
\end{equation*}%
where $E_{n}$ are familiar Euler numbers.
\end{corollary}

\end{document}